\documentclass[12pt]{article}

      \usepackage{latexsym}
         \usepackage[reqno, namelimits, sumlimits]{amsmath}
         \usepackage{amssymb, amsfonts}
         \usepackage{amsthm}

 \newtheorem{theorem}{Theorem}[section]
 
 \newtheorem{definition}[theorem]{Definition}
 
 \newtheorem{lemma}[theorem]{Lemma}
 
 \newtheorem{remark}[theorem]{Remark}

 \newtheorem{pro}[theorem]{Proposition}
 
\title{Ancient Solutions to Navier-Stokes Equations in Half Space
}

 \author{T. Barker, G. Seregin \\
 Oxford University\\}

\begin{document}
\maketitle

\setcounter{equation}{0}
\section{Introduction}

The goal of the paper is to understand properties of the so-called ancient (backward) solutions to the Navier-Stokes equations. The importance of them in the regularity theory for the Navier-Stokes equations, see, for example, papers \cite {ESS2003}, \cite{KochNadSS2009}, \cite{SerSve2009}, \cite{Ser2012}, and \cite{SerSve2013}, and more generally  in the theory of PDEs is well understood. They  appear as a limit, resulting from rescaling solutions to initial boundary value problems around possible singularities. For the Navier-Stokes equations, this procedure has been described in the above papers.


The weakest version of ancient solutions to the Navier-Stokes equations is as follows.  A vector-valued function $u\in L_{2,loc}(Q_-)$, where $Q_-=\mathbb R^3\times ]-\infty,0[$, is an \textit{ ancient solution} of the Navier-Stokes equations  in $Q_-=\mathbb R^3\times ]-\infty,0[$ if it
 satisfies these equations in the sense of distributions with divergence free test functions, i.e.,
\begin{equation}\label{momentum}
    \int\limits_{Q_-}\Big(u\cdot(\partial_t\varphi+\Delta \varphi)+u\otimes u:\nabla\varphi\Big)dz=0
\end{equation}
for any $\varphi\in C^\infty_{0,0}(Q_-):=\{\varphi\in C^\infty_{0}(Q_-):\,\,\mbox{div}\,\varphi=0\}$ and
\begin{equation}\label{incompressibility}
    \int\limits_{Q_-}u\cdot\nabla q dz=0
\end{equation}
for any $q\in C^\infty_0(Q_-)$.

This class of ancient solutions seems to be too wide. Having in mind the problem of regularity for solutions to the Navier-Stokes equations mentioned above, we can put some additional restrictions in the definition of ancient solutions. 

If an ancient solution $u$ is bounded, we call it a {\it bounded ancient} one. We can go further, see \cite{KochNadSS2009}, and consider an even more narrow class of ancient solutions. We say that a bounded function $u$ is a {\it mild bounded ancient} solution if  $u$
has the following property: for any $A<0$ 
and for $(x,t)\in Q_{A}:= \mathbb{R}^3\times ]A,0[$,
$$u_i(x,t)=\int\limits_{\mathbb R^3}\Gamma(x-y,t-A)u_i(y,A)dy+$$
\begin{equation}\label{mildinwholespace}
+\int\limits^t_A \int\limits_{\mathbb R^3} K_{ijm}(x,y,t-\tau)u_j(y,\tau)u_m(y,\tau)dy d\tau,
\end{equation}
where $\Gamma$ is the  known heat kernel and $K$ is obtained from the Oseen tensor in the following way.  Consider the following boundary value problems
\begin{equation}\label{alpotential}
  \Delta \Phi(x,t)=
  \Gamma(x,t).
\end{equation}
Using $\Phi$, we define
$$K_{mjs}(x,y,t)=\delta_{mj}\frac{\partial^3 \Phi}{\partial y_i\partial y_i\partial y_s}(x,y,t)-\frac{\partial^3 \Phi}{\partial y_m\partial y_j\partial y_s}(x,y,t),$$ where $\delta_{mn}$ is Kronecker's symbol.

It has been shown in \cite{KochNadSS2009}, that any mild bounded ancient solution  is infinitely smooth in space-time.

One can give an equivalent definition of mild bounded ancient solutions.
\begin{pro}\label{equivalence in space} A bounded function $u$ in $Q_-$ is a mild bounded ancient solution to the Navier-Stokes equation if and only if there is a pressure $p\in L_\infty(-\infty,0;BMO)$ such that the pair $u$ and $p$ satisfy the Navier-Stokes equations in the sense of distributions. \end{pro}
This statement seems to be known and we give its prove for completeness. 

One of the interesting consequences of the above proposition is an alternative proof of smoothness of mild bounded ancient solutions, see \cite{Ser2014}.

The conjecture that has been made in \cite{KochNadSS2009}
reads: any mild bounded ancient solution is a constant in
$Q_-$. The validity of this conjecture is known in several cases, see details in \cite{KochNadSS2009} and \cite{SerSve2009}. The connection with a possible blowup of a solution to the initial  value problem
$$\partial_tv+v\cdot\nabla v-\Delta v=-\nabla q,\qquad {\rm div}\, v=0$$
in $Q_\infty=\mathbb R^3\times ]0,\infty[$,
$$v(\cdot,t)=u_0(\cdot)\in C^\infty_{0,0}(\mathbb R^3_+)=\{v\in C^\infty_0(\mathbb R^3):\,\,{\rm div}\,v=0\},$$
is as follows. Assume that there is a blowup at $t=T$, i.e.,
$$\|v(\cdot,t)\|_{\infty,\mathbb R^3}\to\infty$$
as $t\to T_{-}$. Then there exists a mild bounded ancient solution $u$ with $|u(0)|=1$. If the aforesaid conjecture is true then $u(x,t)=c$, where $c$ is a constant vector such that $|c|=1$. This would rule out blowups of Type I for which a certain scale-invariant quantity is bounded.

Now, let us formulate the main results of the paper about mild bounded ancient solutions to the Navier-Stokes equations in half space, starting with a definition of  distributional ancient solutions. From now on we denote $Q^{+}_{-}:= \mathbb{R}^{3}_{+}\times ]-\infty,0[$. 
We say that $u\in L_2(B_+(R))$ for any $R>0$ is an ancient solution  if $u$
satisfies
\begin{equation}\label{distrbutional in half}
    \int\limits_{Q^+_-}\Big(u\cdot(\partial_t\varphi+\Delta \varphi)+u\otimes u:\nabla\varphi)dx dt=0
\end{equation}
for any $\varphi\in C^\infty_{0,0}(Q_-)$ with $\varphi(x',0,t)=0$ for any $x'\in \mathbb R^2$ and for any $t<0$.
 Moreover $u$ satisfies 
\begin{equation}\label{incompressibility in hs}
    \int\limits_{Q_-^+}u\cdot\nabla q dz=0
\end{equation}
for any $q\in C^\infty_0(Q_-)$.

We can notice that (\ref{distrbutional in half}) and (\ref{incompressibility in hs}) is a weak form for the following
$$\partial_tu+u\cdot \nabla u-\Delta u=-\nabla p,\qquad {\rm div}\, u=0$$
in $Q^{+}_{-}$ for some distribution $p$,
$$u(x',0,t)=0$$
for any $x'\in\mathbb R^2$ and any $-\infty<t<0$.

We shall say that $u$ is a bounded ancient solution in half space if it is  ancient and bounded.

From now on define $Q_{A}^+:= \mathbb{R}^3_{+}\times]A,0[.$
In order to proceed further, we need to recall how one can construct a solution
to the following boundary value for the Stokes equations in half space:
$$\partial_t v-\Delta v +\nabla q=f,\qquad {\rm div }\,v=0$$
in $Q_{A}^{+}$,
$$v(x',0,t)=0$$
for all $x'\in \mathbb R^2$ and $t\in ]A,0[$,
$$v(\cdot,A)=u_0(\cdot)$$
in $\mathbb R^3_+$.
It is  assumed that $f$ and $u_0$ are divergence free and $f_3(x',0,t)=0$. Then a formal solution to the above initial boundary value problem is:
$$v_{i}(x,t)=\int\limits_{\mathbb R_+^3}G_{ij}(x,y,t-A)u_{0j}(y)dy+\int \limits_A^t\int \limits_{\mathbb R^3_+}G_{ij}(x,y,t-s)f_{j}(y,s)dyds.$$
 The Green's function $G$ has been derived by Solonnikov in \cite{Sol1973} and is as follows
\begin{equation}\label{GreenFunction}
    G=G^1+G^2,
\end{equation}
where
$$G^1_{ij}(x,y,t)=\delta_{ij}\Big(\Gamma(x-y,t)-\Gamma(x-y^*,t)\Big),$$
$$G^2_{i\beta}(x,y,t)=
4
\frac {\partial}{\partial x_\beta}\int\limits_0^{x_3}\int\limits_{\mathbb R^2}\frac{\partial E}{\partial x_i}(x-z)\Gamma(z-y^*,t)dz,\quad G^2_{i3}(x,y,t)=0,$$
$y^*=(y',-y_3)$, and $E(x)$ is fundamental solution to the Laplace equation in $\mathbb R^3$.

Let us introduce another potential $K=(K_{mjs})$,
$$K_{mjs}(x,y,t)=\frac{\partial^3 \Phi_{mj}}{\partial y_i\partial y_i\partial y_s}(x,y,t)-\frac{\partial^3 \Phi_{mn}}{\partial y_n\partial y_j\partial y_s}(x,y,t),$$
where the tensor $\Phi=(\Phi_{ij})$ are defined as  solutions to the following boundary value problems
\begin{equation}\label{alpotentialhs}
  \Delta_y \Phi_{mn}(x,y,t)=G_{mn}(x,y,t)
\end{equation}
with $\partial{\Phi_{mn}}/\partial{y_3}(x,y,t)=0$ if $n<3$ and with $\Phi_{mn}(x,y,t)=0$ if $n=3$ at $y_3=0$.


Now, we are in position to  define mild bounded ancient solutions in a half space.

\begin{definition}\label{defmildbddansol in half} A bounded divergence free function $u$ in $Q_-^+$ is called a mild bounded ancient solution if,  for any $A<0$ 
and
any $(x,t)\in\mathbb Q^+_A$,
$$u_i(x,t)=\int\limits_{\mathbb R^3_+}G_{ij}(x-y,t-A)u_i(y,A)dy+$$
\begin{equation}\label{mildinhalfspace}
+\int\limits^t_A \int\limits_{\mathbb R^3_+} K_{ijm}(x,y,t-\tau)u_j(y,\tau)u_m(y,\tau)dy d\tau.
\end{equation}
  \end{definition}

  To state our main result, we need to introduce the following operator.
  Given $H=(H_{ij})\in L_\infty(\mathbb R^3_+)$, there exists a unique function $p^1\in L_2(B_+(R))$ for any $R>0$ with $[p^1]_{B_+}=0$ with the following properties: the even extension of it to $\mathbb R^3$ belongs to the space $BMO$,
$$\int\limits_{\mathbb R^3_+}p^1\Delta \varphi dx=-\int\limits_{\mathbb R^3_+}H:\nabla^2\varphi dx$$
for any  $\varphi \in C^\infty_0(\mathbb R^3_+)$ with $\varphi_{,3}(x',0)=0$ and
$$\|p^1\|_{BMO}\leq A\|H\|_{\infty,\mathbb R^3_+},$$
where $A$ is an absolute constant. We shall use notation $p^1:=p^1_H$. 

We notice that
 if $H=(H_{ij})$ is sufficiently smooth 
 and vanishes on the boundary $x_3=0$, the function $p^1_H$ is a solution to the  Neumann boundary value problem:
$$\Delta p^1_H=-{\rm div}{\rm div}H $$
in $\mathbb R^3_+$ and
$$p^1_{H,3}(x',0)=0.$$

\begin{theorem}\label{mabssmooth}
Suppose $u\in L_{\infty}(Q_{-}^{+})$ is an arbitrary  mild bounded ancient solution in $Q_-^+$.
Then $u$ is of class $C^{\infty}$ and moreover
$$\sup_{(x,t)\in Q_{-}^{+}}(|\partial_{t}^{k}\nabla^{l}u(x,t)|+|\partial_{t}^{k}\nabla^{l+1}p(x,t)|)+$$$$+\|\partial_{t}^{k}p^{1}\|_{L_{\infty}(BMO)}\leqslant C(k,l,\|u\|_{L_{\infty}(Q_{-}^{+})})<\infty$$
for any k, l\,$\mathrm{= 0,1}$\ldots. Here, $p^{1}=p^1_{u\otimes u}$.
\end{theorem}


\begin{theorem}\label{hsmildsol}  A bounded function $u$ is a mild bounded ancient solution if and only if  there exists
a pressure $p$ such that $p=p^1_{u\otimes u}+p^2$,
where  $p^2(\cdot,t)$ is a harmonic function in $\mathbb R^3_+$ whose gradient satisfies the estimate
\begin{equation}\label{log}
  |\nabla p^2(x,t)|\leq c\ln(2+ 1/{x_3})
\end{equation}
for all $(x,t)\in Q_-^+$.
Morevoer, $p^2$ has the property
 \begin{equation}\label{p2}
   \sup\limits_{x'\in\mathbb R^2}|\nabla p^2(x,t)|\to 0
 \end{equation}
as $x_3\to \infty$ and for any $t<0$;  $u$ and $p$ satisfy  (\ref{incompressibility in hs}) and
\begin{equation}\label{momentumpres}
    \int\limits_{Q^+_-}\Big(u\cdot(\partial_t\varphi+\Delta \varphi)+u\otimes u:\nabla\varphi+p{\rm div}\,\varphi\Big)dx dt=0
\end{equation}
for any $\varphi\in C^\infty_{0}(Q_-)$ with $\varphi(x',0,t)=0$ for any $x'\in \mathbb R^2$ and for any $t<0$.
\end{theorem}

In \cite{SerSve2013}, there has been conjectured that any mild bounded ancient solution is identically equal to zero in $Q_-^+$. At the moment of writing the paper, there are two cases in which the above conjecture is true, see \cite{GHM2014} and \cite{Ser2014-2}. Both cases are two-dimensional and additional scale-invariant assumptions have been imposed. In the first paper vorticity preserves its sign, while in the second one kinetic energy is bounded.

Now, let us consider the following initial boundary value problem
$$\partial_tv+v\cdot\nabla v-\Delta v=-\nabla q,\qquad {\rm div} v=0$$
in $Q^+_\infty=\mathbb R_+^3\times ]0,\infty[$,
$$v(x',0,t)=0$$
for any $x'\in\mathbb R^2$ and $t\in [0,\infty[$, and
$$v(\cdot,t)=u_0(\cdot)\in C^\infty_{0,0}(\mathbb R^3_+)=\{v\in C^\infty_0(\mathbb R^3_+):\,\,{\rm div}\,v=0\}$$
 Suppose that there is a blowup at $t=T$, i.e.,
$$\|v(\cdot,t)\|_{\infty,\mathbb R^3_+}\to\infty$$
as $t\to T_{-}$. Then there exists a sequence $z_n=(x^{(n)},t_n)$ such that $t_n>0$, $t_n\to T_{-}$, and $$M_n=|u(z^{(n)})|=\sup\limits_{0<t\leq
t_n}\sup\limits_{x\in\mathbb R^3_+}|u(x,t)|\to\infty.$$
If $x_3^{(n)}M_n\to\infty$, there exists a mild bounded ancient solution $u$ in the whole space such that $|u(0)|=1$. If $x_3^{(n)}M_n\to a<\infty$, there exists a mild bounded ancient solution in a half space such that $|u(a)|=1$.

In conclusion, we notice that the validity of both conjectures allow us to rule out at least Type I blowups of solution to initial boundary value problem for the Navier-Stokes equations in half space.

\setcounter{equation}{0}
\section{Proof of Theorem   \ref{mabssmooth}}

Before starting the proof, we remind  known facts (due to Solonnikov, see \cite{Sol1973} and \cite{Sol2003}), about the Green function and the kernel $K$.

It is not so difficult to see that the kernel $K$ has the 
structure
\begin{equation}\label{kernelK}
K_{ism}(x,z,t)=\overline{K}_{ism}(x,z,t)+\widehat{K}_{ism}(x,z,t),
\end{equation}
where $\overline{K}_{ism}(x,z,t)$ is a linear combination of the terms
$$\frac{\partial G_{ij}}{\partial z_k}(x,z,t)$$ and $\widehat{K}_{ism}(x,z,t)$ is a linear combination of the terms
$$\frac {\partial^2}{\partial x_\alpha\partial x_\beta}\int\limits_{\mathbb R^3_+}G_{ij}(x,y,t)\frac {\partial N^{(\pm)}}{\partial y_s}(y,z)dy.$$
Here, $N^{(\pm)}(x,y)=E(x-y)\pm E(x-y^*)$. 

The following estimates for $G^i$  and $\widehat{K}$ have been obtained in different papers  of Solonnikov, see \cite{Sol1973} and \cite{Sol2003}:
$$ \Big|\frac {\partial^{|\alpha|+|\gamma|}G^2}{\partial x^\alpha\partial y^\gamma}(x,y,t-A)\Big|\leq c(\alpha,\gamma) (t-A)^{-\frac {\gamma_3}2}(t-A+x^2_3)^{-\frac{\alpha_3}2}\times $$\begin{equation}\label{GreenFunctEst2}
 \times  (|x-y^*|^2+t-A)^{-\frac {3+|\alpha'|+|\gamma'|}2}\exp{\Big(-\frac {cy^2_3}{t-A}\Big)},
\end{equation}
where $\alpha'=(\alpha_1,\alpha_2)$, $\gamma'=(\gamma_1,\gamma_2)$, and $|\gamma|=0$ or  $1$,
\begin{equation}\label{dergreen2}
    \Big|\frac{\partial G^1_{ij}}{\partial y_i}(x,y,t) \Big|+ |\widehat{K}_{ism}(x,y,t)|\leq \frac c{(|x-y|^2+t)^2}.
\end{equation}
\begin{equation}\label{derivativeintime}\Big|\partial^l_t G^2(x,y,t)\Big|\leq\frac c {t^l(|x'-y'|^2+x_3^2+y^2_3 +t)^\frac 32}\exp{\Big(-\frac {cy^2_3}{t}\Big)}. \end{equation}
for $l=0$ or  $1$.

Let  $K^1$ and $K^2$ be generated by $G^1$ and $G^2$, respectively.  In particular, we have the estimate
\begin{equation}\label{hatK2}
  |\widehat{K}^2(x,y,t)|\leq \frac c{(|x-y^*|^2+t)^2}.
\end{equation}

In what follows, we are going to use special approximations:
\begin{equation}\label{approximationu2}
u_{(k)}(x,t):= \int\limits_{A-1}^{0}\eta_{\frac{1}{k}}(t-\tau)\int\limits_{\mathbb{R}^{3}_{+}}\omega_ {\frac{1}{k}}(x-y)\phi_{k}(y)u^{(\frac{2}{k})}_{}(y,\tau)dyd\tau
\end{equation}
 and
 \begin{equation}\label{approximationu3}
 u_{(k)A}(x):= \int\limits_{\mathbb{R}^{3}_{+}}\omega_{\frac{1}{k}}(x-y)u^{(\frac{2}{k})}(y,A)dy.
 \end{equation}
Here, $u^{(h)}_{i}(y,s)=u_{i}(y',y_3-h,s)$ if $y_3>h$ and $u^{(h)}_{i}(y,s)=0$ if $0<y_3\leq h$. The function $\phi_{k}\in C_{0}^{\infty}(B(k+1))$, $\phi_{k}\equiv 1$ on $B(k)$, has  the additional property that the bounds of $D^{\alpha}\phi_{k}$ only depend on $|\alpha|$. The standard mollifiers are denoted by $\eta$ and $\omega$, repectively. The properties of the approximation scheme  are that $u_{(k)}\in C_{0}^{\infty}(]A-2,1[\times\mathbb{R}_{+}^{3})$ and that (up to subsequence) $$u_{(k)}\otimes u_{(k)}\stackrel{*}{\rightharpoonup} u\otimes u $$
in $L_\infty(Q^+_A;\mathbb{M}^{3\times 3})$. It is noticed that $u_{(k)A}$ is a smooth solenodial vector field, with bounded derivatives that all vanish near $x_{3}=0$. Furthermore,
$$u_{(k)A}\stackrel{*}{\rightharpoonup} u(y,A) $$
in $L_{\infty}(\mathbb R^3_{+};\mathbb{R}^{3})$.

We let  $F^{k}:= u_{(k)}\otimes u_{(k)}$ 
and then
$$U_{(k)}(x,t):= \int\limits_{\mathbb{R}^{3}_{+}}G(x,y,t-A)u_{(k)A}(y) dy+\int\limits_{A}^{t}\int\limits_{\mathbb{R}^{3}_{+}}K(x,y,t-\tau)F^k(y,\tau) dyd\tau.
$$
It is not so difficult to infer that (up to subsequence):
\begin{equation}\label{weak*convergapprox}U_{(k)}\stackrel{*}{\rightharpoonup} u
\end{equation}
in $L_{\infty}(Q_{A}^{+};\mathbb{R}^{3})$.

 To treat the second term on the right hand side of representation formula, we are going to use the following statement.

\begin{lemma}\label{itegbyparts} Let $ F\in W^1_\infty(\mathbb R^3_+)\cap C^{1}(\overline{B}_+(R))$ for any $R>0$ with $F=0$ and $F_{3j,j}=0$ on the plane $x_3=0$. In addition, assume that ${\rm div}{\rm div}F\in L_\infty(\mathbb R^3_+)\cap C(\overline{B}_+(R))$ for any $R>0 $. Then the  identity
$$\int\limits_{\mathbb R^3_+}\Delta_y\Phi_{ij}(x,y,t)f_j(y)dy =\int\limits_{\mathbb R^3_+}G_{ij}(x,y,t)f_j(y)dy =$$
$$=\int\limits_{\mathbb R^3_+}K_{ijm}(x,y,t)F_{jm}(y)dy.$$
is valid. Here, $f=-{\rm div}F-\nabla p^1_F$.
\end{lemma}
\begin{remark}\label{gradofp1}
Under assumptions imposed on tensor-valued function $F$,
$$\nabla p^1_F\in L_\infty(\mathbb R^3_+)\cap C(\overline{B}_+(R))$$for any $R>0 $.\end{remark}

\textsl{Proof of Lemma \ref{itegbyparts}}. Obviously, we can find a sequence $F^m\in C^\infty_0(\mathbb R^3_+)$ such that
$$F^m,\,\,\nabla F^m,\,\,{\rm div}{\rm div}F^m\stackrel{*}{\rightharpoonup} F,\,\,\nabla F,\,\,{\rm div}{\rm div}F$$
in $L_\infty(Q^+_-)$, respectively. In order to construct such a sequence, we proceed as follows. Let $F^{(h)}(x)=F(x',x_3-h)$ if $x_3>h$ and $F^{(h)}(x)=0$ if $0<x_3\leq h$. Then we let $F^{(h,R)}(x)=\varphi_R(x)F^{(h)}(x)$ with a stadart cut-off function $\varphi_R(x)=\varphi(x/R)$, where $\varphi \in C^\infty_0(B(2))$ and $\varphi\equiv 1$ in $B$. And finally we can produce $F^m$ using $(F^{(h,R)})_\varrho$ with $0<\varrho <h$, where $(g)_\varrho$ is a  mollification of the function $g$.

We also can state that $p^1_{F^m}$ has   decay $\frac{1}{|x|^2}$ as $|x|\to \infty$. So, we do not need to take care of integrability of functions involved because of Solonnikov estimates and the decay of the pressure. Similarly the decay of the pressure allows one to rigorously justify the integration by parts shown below. This is done by proving slow decay of the kernels, using arguments in \cite{Sol2003} and \cite{Sol2003UMN}.

Now, letting $f^m=-{\rm div}F^m-\nabla p^1_{F^m}$, we have
$$A_i:=\int\limits_{\mathbb R^3_+}\Delta_y\Phi_{ij}(x,y,t)f^m_j(y)dy =\int\limits_{\mathbb R^3_+}G_{ij}(x,y,t)f^m_j(y)dy =$$
$$=\int\limits_{\mathbb R^3_+}\Phi_{ij,kk}(x,y,t)(-F^m_{js,s}(y) -p^1_{F^m,j}(y))dy=$$
$$=\int\limits_{\mathbb R^3_+}\Phi_{ij,kks}(x,y,t)F^m_{js}(y)dy +\int\limits_{\mathbb R^2}\Phi_{ij,3}(y',0)p^1_{F^m,j}(y',0)dy'+$$
$$+\int\limits_{\mathbb R^3_+}\Phi_{ij,k}(x,y,t)p^1_{F^m,jk}(y))dy.$$
By our assumptions on boundary values of functions $\Phi$ and $p^1_{F^m}$ and their derivatives, the integral over the plane $x_3=0$ vanishes. So, we have
$$A_i=\int\limits_{\mathbb R^3_+}\Phi_{ij,kks}(x,y,t)F^m_{js}(y)dy-\int\limits_{\mathbb R^2}\Phi_{ij}(y',0)p^1_{F^m,3j}(y',0)dy'-$$
$$-\int\limits_{\mathbb R^3_+}\Phi_{ij}(x,y,t)\Delta p^1_{F^m,j}(y))dy. $$
For the same reason, the surface integral is equal to zero and using the pressure equation, we find
$$A_i=\int\limits_{\mathbb R^3_+}\Phi_{ij,kks}(x,y,t)F^m_{js}(y)dy+\int\limits_{\mathbb R^3_+}\Phi_{ij}(x,y,t)F^m_{sk,jsk}(y)dy=$$
$$=\int\limits_{\mathbb R^3_+}\Phi_{ij,kks}(x,y,t)F^m_{js}(y)dy-\int\limits_{\mathbb R^3_+}\Phi_{ij,jsk}(x,y,t)F^m_{sk}(y)dy. $$

So, the formula of the lemma proved for $F^m$, i.e., we have
$$\int\limits_{\mathbb R^3_+}G_{ij}(x,y,t)f^m_j(y)dy=\int\limits_{\mathbb R^3_+}K_{ijk}(x,y,t)F^m_{jk}(y)dy.$$
Now, the identity of the lemma can be obtained by passage to the limits in the latter identity  as  $m\to \infty$. $\Box$

So, if we let $p^{1(k)}:= p^{1}_{u_{(k)}\otimes u_{(k)}} $, when we have
\begin{equation}\label{approxformula}
U_{i(k)}(x,t)=\int\limits_{\mathbb{R}^{3}_{+}}G_{ij}(x,y,t-A)u_{j(k)A}(y) dyS-$$$$-\int\limits^t_{A}\int\limits_{\mathbb{R}^{3}_{+}}G_{ij}(x,y,t-\tau)\Big[\frac{\partial}{\partial y_{l}}F^{k}_{jl}(y,\tau)+\frac{\partial}{\partial y_{j}}p^{1(k)}(y,\tau)\big]dyd\tau.
\end{equation}

We then put
$U_{(k)}=U_{(k)}^{1}+U_{(k)}^{2}$ according to splitting of the kernel described in (\ref{GreenFunction}). Furthermore, decompose $U_{(k)}^{1}= U_{(k)}^{1,1}+U_{(k)}^{1,2}$ and $U_{(k)}^{2}= U_{(k)}^{2,1}+U_{(k)}^{2,2}$.
Where for $m=1,2$:
\begin{equation}\label{approxsplit1}
U_{i(k)}^{m,1}(x,t):= \int\limits_{\mathbb{R}^{3}_{+}}G_{ij}^{m}(x,y,t-A)u_{j(k)A}(y) dy
\end{equation}
and
\begin{equation}\label{approxsplit2}
U_{i(k)}^{m,2}(x,t):= -\int\limits^t_{A}\int\limits_{\mathbb{R}^{3}_{+}}G_{ij}^{m}(x,y,t-\tau)\Big[\frac{\partial}{\partial y_{l}}F^{k}_{jl}(y,\tau)+\frac{\partial}{\partial{y_j}}p^{1(k)}(y,\tau)\big]dyd\tau.
\end{equation}

Notice, that we may integrate by parts in (\ref{approxsplit2}) to get:
$$U_{i(k)}^{m,2}(x,t):= \int\limits^t_{A}\int\limits_{\mathbb{R}^{3}_{+}}\frac{\partial}{\partial y_{l}}G_{ij}^{m}(x,y,t-\tau)\Big[F^{k}_{jl}(y,\tau)+\delta_{jl}p^{1(k)}(y,\tau)\Big]dyd\tau.$$
This is permissible by the facts that $G^{2}_{i3}(x,y,t)=0$, $G^{1}_{ij}(x,y,t)=0$ (on $y_{3}=0$), the spatial decay of $G_{1}$ and $G_{2}$ and that the approximation scheme implies that $p^{1(k)}(y,\tau)$ has spatial decay of order $|y|^{-2}$.
\par Now we proceed in proving the main body of Theorem \ref{mabssmooth}. Most of the proof is split into four main Propositions. The first two Propositions are derived from arguments from \cite{SerSve2013}. However, we provide some adjustments, simplifications and demonstrate how those arguments interact with the aforementioned approximation scheme.

In what follows, we are going to use additional notation.
For $p$ and $q$ between 1 and infinity we say that $f\in L_{p,q,unif}(Q_{A}^{+})$, 
 if
$$\|f\|_{L_{p,q,unif}(Q_{A}^{+})}^{q}:= \sup_{x\in\mathbb{R}_{+}^{3}}\int\limits^0_A\Big(\int\limits_{B^{+}(x,1)}|f(y,\tau)|^{p}dy\Big)^{\frac{q}{p}}d\tau< \infty,$$
where $B_+(x,1):=\{y\in B(x,1): \,y_3>x_3\}.$ In addition, for $-\infty<C<D<\infty$, we will denote: $Q_{C,D}^{+}:=\mathbb{R}^3_+\times ]C,D[.$
From now on we use the terms even and odd extensions to mean the following. For $f:\mathbb{R}^{3}_+\rightarrow\mathbb{R}$, define $f_{even}:\mathbb{R}^3\rightarrow\mathbb{R}$ by $f_{even}(y):= f(y)$ for $y_3\geqslant 0$ and $f_{even}(y):= f(y^*)$ for $y_{3}<0$. This is referred to as the even extension of $f$.
The odd extension of $f$ is similarly defined.

\begin{pro}\label{gradbdd}
Suppose $u\in L_{\infty}(Q_{-}^{+})$ satisfies all the hypothesis of Theorem \ref{mabssmooth}.
Then the following is satisfied:
$$\sup_{(x,t)\in Q_{-}^{+}}(|\nabla u(x,t)|+|\nabla p^{1}(x,t)|)
+\|p^{1}\|_{L_{\infty}(BMO)})\leqslant C(\|u\|_{L_{\infty}(Q_{-}^{+})})<\infty$$

\end{pro}
\textsl{Proof of Proposition \ref{gradbdd}}\ For brevity let $d:=\|u\|_{L_{\infty}(Q_{A}^{+})}$.
Notice that by classical singular integral theory, we get for the even extension of the pressure:
\begin{equation}\label{pressureBMO}
\sup_{k}\|p^{1(k)}\|_{L
_{\infty}(]A,0[;BMO(\mathbb{R}^{3})}\leqslant C(d).
\end{equation}
 By the proerties of the  heat kernel and estimates (\ref{GreenFunctEst2}) it is obtained that
\begin{equation}\label{stokepotentbdd}
\sup_{k,\, (x,t)\in Q_{\frac{A}{2}}^{+}}|\nabla^{|\alpha|}U_{(k)}^{1,1}(x,t)|+|\nabla^{|\alpha|}U_{(k)}^{2,1}(x,t)|\leqslant C(d,A,|\alpha|).
\end{equation}
From (\ref{pressureBMO}), along with arguments in \cite{SerSve2013} (see Lemma 6.1 there), obtain also that (for $(x,t)\in Q_{A}^{+}$):

\begin{equation}\label{solonnikovpartest}
|\nabla^{|\alpha|}U_{(k)}^{2,2}(x,t)|\leqslant C(|\alpha|,d)\int\limits^{\sqrt{-A}}_{0}\frac{dq}{(x_{3}^{2}+q^{2})^{\frac{|\alpha|}{2}}}.
\end{equation}
Hence, $\|\nabla U_{(k)}^{2}\|_{L_{s,unif}(Q_{\frac{A}{2}}^{+})}+\|\nabla U_{(k)}^{1,1}\|_{L_{s,unif}(Q_{\frac{A}{2}}^{+})}\leqslant C(d,s,A)$. By the definition of $G^1$, the following is satisfied in $\mathcal{D}^{'}(Q_{A})$ :
\begin{equation}\label{halfspaceheatpart}
\partial_{t}U^{1,2}_{(k)_{odd}}-\Delta U^{1,2}_{(k)_{odd}}= -{\rm div}\widetilde H,
\end{equation}
where $\widetilde H_{i\alpha}=H^{odd}_{i\alpha}$, $\alpha=1,2$, and $\widetilde H_{i3}=H^{even}_{i3}$,
 $i=1,2,3$,
 $$H^{k}_{ij}(x,t):= F^{k}_{ij}(x,t)+ \delta_{ij}(p^{1(k)}-[p^{1(k)}]_{B(z,2)}(t)).$$
Here, $z\in\mathbb{R}^{3}_{+}$ is arbitrary and $[p^{1(k)}]_{B(z,2)}(t)$ signifies the average over $B(z,2)$ for the even extension.  Using (\ref{pressureBMO}) and local regularity theory for heat equation (e.g Appendix of \cite{NRS1996}), obtain:
$$\|\nabla U_{(k)_{odd}}^{1,2}\|_{L_{s}(B(z,1)\times]\frac{A}{2},0[)}\leqslant C(d,s,A).$$
Thus: $\|\nabla U_{(k)}\|_{L_{s,unif}(Q_{\frac{A}{2}}^{+})}\leqslant C(d,s,A).$ Notice that for $x'\in\mathbb{R}^2$, $y\in\mathbb{R}^3_+$ and $t>0$, we have
$$\lim_{\varepsilon\rightarrow 0^+}G^{2}_{ij}(x',\epsilon,y,t)=0.$$ Thus, using additional properties of the heat kernel and the properties of the approximations $u_{(k)}$ and $p^{1(k)}$, it can be obtained that $\nabla U_{i(k)}(x,t)$ is bounded in $\mathbb{R}^{3}_{+}\times ]A,0[$ . Furthermore, for $(x',t)\in\mathbb{R}^{2}\times]A,0[$ one obtains
$$\lim_{\varepsilon\rightarrow 0}U_{i(k)}(x',\epsilon,t)=0.$$
Hence, we obtain a weak derivative formula for classes of test functions not necessarily zero on $x_{3}=0$. That is, for $\varphi\in C_{0}^{\infty}(Q_{A})
$:
\begin{equation}\label{zerotraceofapprox}
\int\limits^0_{A}\int\limits_{\mathbb{R}^{3}_{+}}U_{i(k)}(y,\tau)\partial_{j}\varphi(y,\tau) dyd\tau= -\int\limits^0_{A}\int\limits_{\mathbb{R}^{3}_{+}}\partial_{j}U_{i(k)}(y,\tau)\varphi(y,\tau) dyd\tau.
\end{equation}
So using (\ref{weak*convergapprox}) one has that  $\|\nabla u\|_{L_{s,unif}(Q_{{A-1}}^{+})}\leqslant C(d,s,A)$ and that (\ref{zerotraceofapprox}) holds also for $u$ (note one can replace $Q_{A}^{+}$ with $Q_{-}^{+}$ and $A$ with $-\infty$ here). Next, fix $\delta<0$, with $|\delta|$ small and let $k$ be sufficiently large such that $\delta+\frac{1}{k}<0$.
Observing the structure of the approximations $u_{(k)}$, the analogue of (\ref{zerotraceofapprox}) (with $u$) gives for
$(x,t)\in Q_{A,\delta}^+$:
\begin{equation}\label{derapprox}
\partial_{j}u_{i(k)}(x,t)= \int\limits_{A-1}^{0}\eta_{\frac{1}{k}}(t-\tau)\int\limits_{\mathbb{R}^{3}_{+}}\omega_ {\frac{1}{k}}(x-y)\frac{\partial}{\partial y_{j}}\phi_{k}(y)u^{(\frac{2}{k})}_{i}(y,\tau)dyd\tau+
\end{equation}
$$+\int\limits_{A-1}^{0}\eta_{\frac{1}{k}}(t-\tau)\int\limits_{\mathbb{R}^{3}_{+}}\omega_ {\frac{1}{k}}(x-y)\phi_{k}(y)(\partial_{j}u_{i})^{(\frac{2}{k})}(y,\tau)dyd\tau.
$$
 Thus the improvement of $u$ gives  (for sufficiently large $k\geqslant K(\delta)$):
\begin{equation}\label{approxgradest}
\|\nabla u_{(k)}\|_{L_{s,unif}(Q_{A,\delta}^+)}\leqslant C(d,s,A).
\end{equation}
Now $p^{1(k)}$ satisfies (for appropriate even and odd extensions of $p^{1(k)}$ and $u_{i(k)}u_{j(k)}$)
$$\Delta p^{1(k)}(x,t)= -{\rm div}{\rm div}( u_{(k)}(x,t)\otimes u_{(k)}(x,t))$$ in $Q_{A}$. Local regularity theory for Laplace equation gives:

\begin{equation}\label{approxpressureest}
\|\nabla p^{1(k)}\|_{L_{s,unif}(Q_{A,\delta}^+)}\leqslant C(d,s,A).
\end{equation}
Using local regularity theory for heat equation, we  find from (\ref{halfspaceheatpart}):
\begin{equation} \label{estheatpart1}
\|\partial_{t}U_{i(k)_{odd}}^{1,2}\|_{L_{s}(B(z,1)\times ]\frac{3A}{8},\delta[)}+\|\nabla^{2}U_{i(k)_{odd}}^{1,2}\|_{L_{s}(B(z,1)\times ]\frac{3A}{8},\delta[)}\leqslant C(d,s,A).
\end{equation}
For $s$ sufficiently large it is seen that $\nabla U_{i(k)}^{1,2}$ is bounded (in fact H\"older continuous) in $Q_{\frac{3A}{8},\delta}^+$ with
\begin{equation}\label{gradheatpartbdd}
\sup_{Q_{\frac{3A}{8},\delta}^+}|\nabla  U_{i(k)}^{1,2}(x,t)|\leqslant C(d,s,A).
\end{equation}

To estimate $\nabla U^{2,2}_{(k)}$ we need the following statement whose proof is  contained in the Appendix.  

 \begin{lemma}\label{Uniforminteggeneralise}
 Suppose $f$ is in $L_{s,l,unif}(Q_{A}^{+})$.
 Furthermore  assume that $1<s\leqslant l\leqslant\infty$ along with:
 \begin{equation}\label{Uniformintegcondition}
 \frac{3}{s}+\frac{2}{l}<1.
 \end{equation}
  Then it follows that for $(x,t)$ in $Q_{A}^{+}$
$$\int\limits^t_A\int\limits_{\mathbb{R}^{3}_{+}}|\nabla_{y}G^{2}(x,y,t-\tau)f(y,\tau)| dyd\tau<\infty.$$
Furthermore,
$$\sup_{(x,t)\in Q_{A}^{+}}\int\limits^t_A\int\limits_{\mathbb{R}^{3}_{+}}|\nabla_{y}G^{2}(x,y,t-\tau)f(y,\tau)| dyd\tau\leqslant C(A,s,l)\|f\|_{L_{s,l,unif}(Q_{-A}^{+})}.$$
 \end{lemma}

 \begin{remark}\label{Luniformxgrad}
 Observing (\ref{GreenFunctEst2}), we see that Lemma \ref{Uniforminteggeneralise} holds if we replace $\nabla_{y}$ with $\nabla_{x}$.
 \end{remark}

Now, using Lemma \ref{Uniforminteggeneralise} (in particular Remark \ref{Luniformxgrad}), (\ref{approxgradest}) and (\ref{approxpressureest}), we  find:
\begin{equation}\label{gradsolonnpartbdd}
\sup_{Q^{+}_{A,\delta}}|\nabla  U_{(k)}^{2,2}(x,t)|\leqslant C(d,s,A).
\end{equation}
Hence, $\sup_{Q^{+}_{\frac{A}{2},\delta}}|\nabla  U_{(k)}(x,t)|\leqslant C(d,s,A).$ The conclusion regarding boundedness of the gradient of $u$ in $Q_{-}^{+}$ is inferred from taking limits and  time-shift arguments.

The statement regarding $\|p^{1}\|_{L_{\infty}(BMO)}$ is deduced from (\ref{pressureBMO}). It remains to prove $\sup_{(x,t)\in Q_{-}^{+}}|\nabla p^{1}(x,t)|\leqslant C(\|u\|_{L_{\infty}(Q_{-}^{+})})$. Notice from (\ref{incompressibility in hs}) that for $(x,t)\in Q^{+}_{A,\delta}$:
\begin{equation}\label{divofapprox}
{\rm div}u_{(k)}(x,t)= \int\limits_{A-1}^{0}\eta_{\frac{1}{k}}(t-\tau)\int\limits_{\mathbb{R}^{3}_{+}}\omega_{\frac{1}{k}}(x-y)\nabla\phi_{k}(y).u^{(\frac{2}{k})}(y,\tau)dyd\tau.
\end{equation}
From the latter, it follows that ${\rm div}\, u_{(k)}$ and $\nabla {\rm div}\, u_{(k)}$ are bounded function in space-time and in $k$. Local regularity for Laplace equation gives (for the even extension of the pressure):
\begin{equation}\label{hesspressureest}
\sup_{t\in  ]A,\delta[}\|\nabla^{2}p^{1(k)}(\cdot,t)\|_{L_{s}(B(z,1))}+\|\nabla p^{1(k)}(\cdot,t)\|_{L_{s}(B(z,1))}\leqslant C(s,d,A).
\end{equation}
The conclusion is then reached by arguments similar to those previously mentioned. Proposition \ref{gradbdd} is proven. $\Box$

Before proceeding the proof of Theorem \ref{mabssmooth}, let us adopt the notation:
 $$\mathbb{R}^{3}_{\gamma}:= \{(x',x_{3})\in\mathbb{R}^{3}: |x_{3}|\geqslant\gamma\},$$
 $$\mathbb{R}^{3}_{\gamma+}:= \mathbb{R}^{3}_{+}\cap \mathbb{R}^{3}_{ \gamma}.$$
\begin{pro}\label{propertytimeder}
Suppose all the assumptions of Proposition \ref{gradbdd} hold.
Then $u$ also satisfies:
\begin{equation}\label{Luniftimeder}
\|\partial_{t} u\|_{L_{s,unif}(Q_{A}^{+})}\leqslant C(d,s,A)
\end{equation} (for any $A\in ]-\infty,0[$).
\begin{equation}\label{bddawaytimeder}
\sup_{x\in\mathbb{R}^{3}_{\gamma+},\, t\in ]-\infty,0[}|\partial_{t} u(x,t)|\leqslant C(d,\gamma).
\end{equation}
\end{pro}
\textsl{Proof of Proposition \ref{propertytimeder}}
First notice that, due to estimates (\ref{GreenFunctEst2}), for $x\in\mathbb{R}_{+}^{3}$:
\begin{equation}\label{Greenfunctionsmalltime}
\int\limits_{\mathbb{R}^{3}_{+}}|G^{2}(x,y,\epsilon)|dy\leqslant C(n)\frac{\sqrt{\epsilon}}{x_{3}}.
\end{equation}

Note that properties of the approximation scheme imply $F^{k}_{jl}$ and $p^{1(k)}$ are smooth on $\mathbb{R}^{3}_{+}\times]A,0[$ with bounded derivatives (the bounds may be dependent on $k$).
Change variables to get:
$$U_{i(k)}^{2,2}(x,t)= \int\limits^{t-A}_{0}\int\limits_{\mathbb{R}^{3}_{+}}G_{ij}^{2}(x,y,\lambda)\Big[\frac{\partial}{\partial y_{l}}F^{k}_{jl}(y,t-\lambda)+\frac{\partial}{\partial{y_j}}p^{1(k)}(y,t-\lambda)\big]dyd\tau.
$$
It follows that differentiation in time is permissible and gives:
$$\frac{\partial}{\partial{t}}U_{i(k)}^{2,2}(x,t)=-\int\limits^{t-A}_{0}\int\limits_{\mathbb{R}^{3}_{+}}G_{ij}^{2}(x,y,\lambda)\frac{\partial}{\partial \lambda}\Big[\frac{\partial}{\partial y_{l}}F^{k}_{jl}(y,t-\lambda)+\frac{\partial}{\partial{y_j}}p^{1(k)}(y,t-\lambda)\big]dyd\tau+$$$$+
\int\limits_{\mathbb{R}^{+}_{3}}G_{ij}^{2}(x,y,t-A)\Big[\frac{\partial}{\partial y_{l}}F^{k}_{jl}(y,A)+\frac{\partial}{\partial{y_j}}p^{1(k)}(y,A)\big]dy.
$$
Notice we have (for $\epsilon>0$):
$$-\int\limits^{t-A}_{\epsilon}\int\limits_{\mathbb{R}^{3}_{+}}G_{ij}^{2}(x,y,\lambda)\frac{\partial}{\partial \lambda}\Big[\frac{\partial}{\partial y_{l}}F^{k}_{jl}(y,t-\lambda)+\frac{\partial}{\partial{y_j}}p^{1(k)}(y,t-\lambda)\big]dyd\lambda=$$$$=-\int\limits_{\mathbb{R}^{+}_{3}}G_{ij}^{2}(x,y,t-A)\Big[\frac{\partial}{\partial y_{l}}F^{k}_{jl}(y,A)+\frac{\partial}{\partial{y_j}}p^{1(k)}(y,A)\big]dy+$$$$+
\int\limits_{\mathbb{R}^{+}_{3}}G_{ij}^{2}(x,y,\epsilon)\Big[\frac{\partial}{\partial y_{l}}F^{k}_{jl}(y,t-\epsilon)+\frac{\partial}{\partial{y_j}}p^{1(k)}(y,t-\epsilon)\big]dy+$$$$+
\int\limits^{t-A}_{\epsilon}\int\limits_{\mathbb{R}^{3}_{+}}\frac{\partial}{\partial\lambda}G_{ij}^{2}(x,y,\lambda)\Big[\frac{\partial}{\partial y_{l}}F^{k}_{jl}(y,t-\lambda)+\frac{\partial}{\partial{y_j}}p^{1(k)}(y,t-\lambda)\big]dyd\lambda.
$$
Thus using (\ref{Greenfunctionsmalltime}) and estimates (\ref{GreenFunctEst2})-(\ref{dergreen2}) obtain:
\begin{equation}\label{solonnparttimeder}
\frac{\partial}{\partial{t}}U_{i(k)}^{2,2}(x,t)= -\int\limits^t_{A}\int\limits_{\mathbb{R}^{3}_{+}}\frac{\partial}{\partial{t}}G_{ij}^{2}(x,y,t-\tau)\Big[\frac{\partial}{\partial y_{l}}F^{k}_{jl}(y,\tau)+\frac{\partial}{\partial{y_j}}p^{1(k)}(y,\tau)\big]dyd\tau.
\end{equation}
We know that  $\|\nabla F^{k}\|_{L_{\infty}(Q^{+}_{A,\delta})}+\|\nabla p^{1(k)}\|_{L_{\infty}(Q^{+}_{A,\delta})}\leqslant C(d)$.
So, using the estimate (\ref{derivativeintime}), it is deduced that for $(x,t)\in Q^{+}_{A,\delta}$:
\begin{equation}\label{logtimederest}
|\partial_{t}U_{(k)}^{2,2}(x',x_{3},t)|\leqslant\log\left(2-A+\frac{(-A)^{2}-A}{x_{3}^{2}}\right)
\end{equation}
and
$$ \sup_{(x,t)\in Q^+_{\frac{A}{2}}}(|\partial_{t}U_{i(k)}^{2,1}(x,t)|+|\partial_{t}U_{i(k)}^{1,1}(x,t)|)\leqslant C(d,A).$$
Thanks to (\ref{estheatpart1}), we  get
\begin{equation}\label{timeLunifapprox}
\|\partial_{t}U_{(k)}\|_{L_{s,unif}(Q^{+}_{\frac{3A}{8},\delta})}\leqslant C(d,s,A).
\end{equation}
In this way (\ref{Luniftimeder}) can be inferred.

To prove (\ref{bddawaytimeder}) it is sufficient to show:
\begin{equation}\label{hessbddawaybndry}
\sup_{x\in\mathbb{R}^{3}_{\gamma+},\,t\in]\frac{A}{4},\delta[}|\nabla^{2} U_{(k)}^{1,2}(x',x_{3},t)|\leqslant c(d,\gamma,A).
\end{equation}
Observing (\ref{stokepotentbdd}), (\ref{solonnikovpartest}) and (\ref{estheatpart1}) it is clear that:
\begin{equation}\label{Lunifhessawaybndry}
\sup_{z\in \mathbb{R}^{3}_{\gamma+}}\|\nabla^{2}U_{(k)}\|_{L_{s}(B(z,\frac{\gamma}{2})\times]\frac{3A}{8},\delta[)}\leqslant c(d,s,A,\gamma).
\end{equation}
By previously mentioned arguments we infer:
$$ \sup_{z\in \mathbb{R}^{3}_{\gamma+}}\|\nabla^{2}u\|_{L_{s}(B(z,\frac{\gamma}{2})\times]A,0[)}\leqslant c(d,s,A,\gamma)$$
and $$ \sup_{z\in \mathbb{R}^{3}_{\gamma+}}\|\nabla^{2}u_{(k)}\|_{L_{s}(B(z,\frac{\gamma}{2})\times]A,\delta[)}\leqslant c(d,s,A,\gamma). $$
Using also (\ref{hesspressureest}), one obtains higher regularity for (\ref{halfspaceheatpart}) through local regularity results for the heat equation. A parabolic imbedding theorem then gives (\ref{hessbddawaybndry}). Proposition \ref{propertytimeder} is proven. $\Box$

The next Proposition  will  improve the previously obtained regularity results. But, first let us state a lemma, which is a simplified version of a more general statement proven in the Appendix, see Lemma \ref{singularintrgal}.
\begin{lemma}\label{singularuniformestsimplified}
 Suppose the measurable kernel $K:\mathbb{R}^{n}\setminus{\{0\}}\rightarrow\mathbb{R}$ satisfies the conditions (see \cite{Stein1970}):
 \begin{equation}\label{convobdd}
 |K(x)|\leqslant B|x|^{-n},\,\,\, {\rm for}\, 0<|x|
 \end{equation}
 \begin{equation}\label{convosmooth}
 \int\limits_{|x|\geqslant 2|y|}|K(x-y)-K(x)|dx\leqslant B,\,\,\, {\rm for}\, 0<|y|
 \end{equation}
 and
 \begin{equation}\label{convocancellation}
 \int\limits_{R_{1}<|x|<R_{2}}K(x) dx= 0,\,\,\, {\rm for}\, 0<R_{1}<R_{2}<\infty.
 \end{equation}

 For suitable $f$ define the singular integral operator :
 \begin{equation} \label{singularintrgalsimplified}
 Tf(x):=\lim_{\epsilon\rightarrow 0}\int\limits_{|x-y|\geqslant\epsilon}K(x-y)f(y)dy.
 \end{equation}

 Take a  compactly supported function $g$ in $L_{p}(\mathbb{R}^{n})$, where $1<p<\infty$. Furthermore assume $g$ is in $L_{\infty}(\mathbb{R}^{n}_{\geqslant 1})$.
 Then it follows that (almost everywhere)
 \begin{equation}\label{decompsimplified}
 Tg(x)= h_{1}(x)+h_{2}(x).
 \end{equation}
 Here,
 \begin{equation}\label{BMOpartsimplified}
 \|h_{1}\|_{BMO(\mathbb{R}^{n})}\leqslant C(n,B)\|g\|_{L_{\infty}(\mathbb{R}^{n}_{\geqslant 1})}
 \end{equation}
 and
 \begin{equation}\label{Luniformpartsimplified}
 \|h_{2}\|_{L_{p,unif}(\mathbb{R}^{n})}\leqslant C(n,p,B)\|g\|_{L_{p,unif}(\mathbb{R}^{n})}.
 \end{equation}

 \end{lemma}
\begin{pro}\label{timederbdd}
Assume $u$ satisfies the all the assumptions of Proposition \ref{gradbdd}.
Then:
\begin{equation}\label{timederpressurebdd}
\sup_{(x,t)\in Q_{-}^{+}}(|\partial_{t}u(x,t)|)+\|\partial_{t}p^{1}\|_{L_{\infty}(BMO)}\leqslant C(\|u\|_{L_{\infty}(Q_{-}^{+})}).
\end{equation}
\end{pro}
\textsl{Proof 
}\
Fix $\delta<0$, with $|\delta|$ small.
From Proposition \ref{propertytimeder}, we have (for $k\geqslant K(\delta$) sufficiently large):
\begin{equation}\label{Luniftimederapprox}
\|\partial_{t}u_{(k)}\|_{L_{s,unif}(Q^+_{A,\delta})}\leqslant C(d,s,A)
\end{equation} and
\begin{equation}\label{timederapproxbddaway}
\sup_{x\in\mathbb{R}^{3}_{\gamma+},\,t\in]A,\delta[}|\partial_{t}u_{(k)}(x,t)|\leqslant C(d,\gamma).
\end{equation}
It is clear that $\partial_{t}p^{1(k)}= p^{1}_{\partial_{t}(u_{(k)}\otimes u_{(k)})}$ and furthermore Lemma \ref{singularuniformestsimplified}  is applicable to a suitable  extension of $\partial_{t}(u_{(k)}\otimes u_{(k)})$ for $n=3$. Hence, it can be written that $\partial_{t}p^{1(k)}= (\partial_{t}p^{1(k)})_{ 1}+(\partial_{t}p^{1(k)})_{2}$.
 Here,
 \begin{equation}\label{BMOpartapproxpres}
 \|(\partial_{t}p^{1(k)})_{ 1}\|_{L_{\infty}(]A,\delta[;BMO(\mathbb{R}^{3})}\leqslant C(d)
 \end{equation}
 and
 \begin{equation}\label{Luniformpartapproxpres}
 \|(\partial_{t}p^{1(k)})_{2}\|_{L_{s,unif}(Q^+_{A,\delta})}\leqslant C(d,s,A).
 \end{equation}
Next,  it  is easy to see that the following is satisfied in $\mathcal{D}^{'}(Q^+_{A,\delta})$ :
\begin{equation}\label{timederivheatodd}
\partial_{t}^{2}U^{1,2}_{(k)_{odd}}-\Delta\partial_{t} U^{1,2}_{(k)_{odd}}= -{\rm div}\widetilde H,
\end{equation}
where $\widetilde H_{i\alpha}=H^{odd}_{i\alpha}$, $\alpha=1,2$, and $\widetilde H_{i3}=H^{even}_{i3}$,
 $i=1,2,3$,
$$H^{k}_{ij}(x,t):= \partial_{t}F^{(k)}_{ij}(x,t)+ \delta_{ij}(\partial_{t}p^{1(k)})_{1}(x,t)-[(\partial_{t}p^{1(k)})_{1}]_{B(\bar{x},2)}(t))+$$$$+ \delta_{ij}(\partial_{t}p^{1(k)})_{2}(x,t).$$
Here $\bar{x}\in\mathbb{R}^{3}_{+}$ is arbitrary.
Hence, (\ref{timeLunifapprox}), (\ref{Luniftimederapprox}), (\ref{BMOpartapproxpres}), (\ref{Luniformpartapproxpres} ), and local regularity of heat equation (for sufficiently large $s>n+2$) give:
\begin{equation}\label{timederhearpartbdd}
\sup_{(x,t)\in Q^+_{\frac{A}{4},\delta}}|\partial_{t}U_{(k)_{odd}}^{1,2}(x,t)|\leqslant C(d,p,A).
\end{equation}

Now, let us examine $\partial_{t}U^{2,2}_{i(k)}$.
One can write:
$$\partial_{t}U_{(k)}^{2,2}(x,t)= -\int\limits_{\mathbb{R}^{3}_{+}}G_{ij}^{2}(x,y,t-A)\Big[\frac{\partial}{\partial y_{l}}F^{(k)}_{jl}(y,A)+\frac{\partial}{\partial{y_j}}p^{1(k)}(y,A)\big]dy+$$$$+ \int\limits_{A}^{t}\int\limits_{\mathbb{R}^{3}_{+}}\frac{\partial}{\partial y_{l}}G_{ij}^{2}(x,y,t-\tau)\Big[\frac{\partial}{\partial t}F^{(k)}_{jl}(y,\tau)+\delta_{jl}\frac{\partial}{\partial t}p^{1(k)}(y,\tau)\Big] dyd\tau. $$
 Proposition \ref{gradbdd} implies that the first term is bounded on $Q^+_{\frac{A}{2}}$ by a constant depending only on $\|u\|_{L_{\infty}(Q_{-}^{+})}$ and $A$. For the second term, decompose as follows:
\begin{equation}\label{timederdecompsolonn}
\int\limits_{A}^{t}\int\limits_{\mathbb{R}^{3}_{+}}\frac{\partial}{\partial y_{l}}G_{ij}^{2}(x,y,t-\tau)\Big[\frac{\partial}{\partial t}F^{(k)}_{jl}(y,\tau)+\delta_{jl}\frac{\partial}{\partial t}p^{1(k)}(y,\tau)\Big] dyd\tau =
\end{equation}
$$= \int\limits_{A}^{t}\int\limits_{\mathbb{R}^{3}_{+}}\frac{\partial}{\partial y_{l}}G_{ij}^{2}(x,y,t-\tau)\Big[\frac{\partial}{\partial t}F^{(k)}_{jl}(y,\tau)+\delta_{jl}\Big(\frac{\partial}{\partial t}p^{1(k)}\Big)_{2}(y,\tau)\Big] dyd\tau+
$$$$+\int\limits_{A}^{t}\int\limits_{\mathbb{R}^{3}_{+}}\frac{\partial}{\partial y_{l}}G_{ij}^{2}(x,y,t-\tau)\delta_{jl}
\Big(\Big(\frac{\partial}{\partial t}p^{1(k)}\Big)_{1}(y,\tau)
-\Big[\Big(\frac{\partial}{\partial t}p^{1(k)}\Big)_{1}\Big]_{B((x^{'},0),a)}(\tau)\Big) dyd\tau.$$
Here, $a= (x_{3}^{2}+t-\tau)^{\frac{1}{2}}$.
For the first part of  (\ref{timederdecompsolonn}) use Lemma \ref{Uniforminteggeneralise} along with estimates  (\ref{Luniftimederapprox}) and (\ref{Luniformpartapproxpres})
to infer that it is bounded on $Q^+_{A,\delta}$ by a constant depending only on $A$ and $d=\|u\|_{L_{\infty}(Q_{-}^{+})}$. For the second part, use (\ref{BMOpartapproxpres}) and arguments from \cite{SerSve2013} to infer that it is bounded on $Q^+_{A,\delta}$ by a constant only depending on $d=\|u\|_{L_{\infty}(Q_{-}^{+})}$.
Thus, putting everything together one has:
\begin{equation}\label{approxtimederbdd}
\sup_{(x,t)\in Q^+_{A,\delta}}|\partial_{t}U_{(k)}(x,t)|\leqslant C(d,A).
\end{equation}
Arguing as before and using  time shift argument, one can get all of the stated conclusions. Proposition \ref{timederbdd} is proven. $\Box$

 The next Proposition briefly describes how the aforementioned arguments can be bootstrapped to obtain analogous statements involving higher time derivatives of $u$.
 \begin{pro}\label{bootstrapintime}
 Suppose $u\in L_{\infty}(Q_{-}^{+})$ satisfies all the assumptions of Theorem \ref{mabssmooth}. Then conclude that:
 \begin{equation}\label{highertimebdd}
 \sup_{(x,t)\in Q_{-}^{+}}(|\partial_{t}^{k}u(x,t)|+|\nabla\partial_{t}^{k} u(x,t)|+|\nabla \partial_{t}^{k}p^{1}(x,t)|)+$$$$+\|\partial_{t}^{k}p^{1}\|_{L_{\infty}(BMO)}\leqslant C(k,l,\|u\|_{L_{\infty}(Q_{-}^{+})})<\infty
\end{equation}
for any $k\mathrm{= 0,1}$\ldots
 \end{pro}
 \textsl{Proof of Proposition \ref{bootstrapintime}}\
We give a brief account of the bootstrap arguments.
Clearly $\partial_{t} u$ also satisfies (\ref{incompressibility in hs}).
By properties of the kernel, if $\varphi\in C_{0}^{\infty}(Q_{A})$ then the following holds:
\begin{equation}\label{timederintegparts}
\int\limits_{A}^{0}\int\limits_{\mathbb{R}^{3}_{+}}\partial_{\tau}^{l}\varphi(y,\tau)U_{i(k)}(y,\tau) dy d\tau= (-1)^{l}\int\limits_{A}^{0}\int\limits_{\mathbb{R}^{3}_{+}}\varphi(y,\tau)\partial_{\tau}^{l}U_{i(k)}(y,\tau) dy d\tau
\end{equation}
and
\begin{equation}\label{timederzerotrace}
\int\limits_{A}^{0}\int\limits_{\mathbb{R}^{3}_{+}}\partial_{y_{q}}\partial_{\tau}^{l}\varphi(y,\tau)U_{i(k)}(y,\tau) dy d\tau= (-1)^{l+1}\int\limits_{A}^{0}\int\limits_{\mathbb{R}^{3}_{+}}\varphi(y,\tau)\partial_{y_{q}}\partial_{\tau}^{l}U_{i(k)}(y,\tau) dy d\tau.
\end{equation}
From Proposition \ref{timederbdd}, we can write:
\begin{equation}\label{timederexpression}
\frac{\partial}{\partial t}U_{i(k)}(x,t)= \frac{\partial}{\partial t}S_i(u_{(k)A})(x,t)
+U_{i(k)}'(x,t,A)
 \end{equation}

$$-\int\limits_{A}^{t}\int\limits_{\mathbb{R}_{+}^{3}}G_{ij}(x,y,t-\tau)\Big[\frac{\partial^{2}}{\partial y_{l}\partial \tau}F_{jl}^{(k)}(y,\tau)+\frac{\partial}{\partial y_{j}}p^{1}_{\frac{\partial}{\partial t}(u_{(k)}\otimes u_{(k)})}(y,\tau)\Big] dyd\tau.
$$
Here, $$U_{i(k)}'(x,t,A)= -\int\limits_{\mathbb{R}^{3}_{+}}G_{ij}(x,y,t-A)\Big[ \frac{\partial}{\partial y_{l}}F^{(k)}_{jl}(y,A)+\frac{\partial}{\partial y_{j}}p^{1(k)}(y,A)\Big] dy$$
and
$$S_{i}(u_{(k)A})(x,t):= \int\limits_{\mathbb{R}^3_{+}}G_{ij}(x,y,t-A)u_{j(k)A}(y)dy.$$
So using Proposition \ref{gradbdd} along with Green function estimates (\ref{GreenFunctEst2}) and (\ref{derivativeintime}) get that:
\begin{equation}\label{surfacetimetermbdd}
\sup_{(x,t)\in Q_{\frac{A}{2}}^{+}}(|\partial_{t}^{k}\nabla^{l}S(u_{A_{k}})(x,t)|+|\partial_{t}^{k}\nabla^{l}U_{i(k)}'(x,t)|)\leqslant C(A,k,l,\|u\|_{L_{\infty}(Q_{-}^{+})})
\end{equation}
any $k,\,l$\,$\mathrm{= 0,1}$\ldots

The third term of (\ref{timederexpression}) is dealt with by splittling the integral according to the kernel decomposition (\ref{GreenFunction}). The arguments are now repeated from Propositions \ref{gradbdd}-\,\ref{timederbdd}. It is possible to repeat this argument indefinitely with higher time derivatives. Proposition \ref{bootstrapintime} is proven. $\Box$
Now, one can recover a pressure $p$ such that in $Q_{-}^{+}$:
$$u\cdot \nabla u-\Delta u+\nabla p=-\partial_tu,\qquad {\rm div}\, u=0$$
in $Q^{+}_{-}$ ,
$$u(x',0,t)=0$$
for any $x'\in\mathbb R^2$ and any $-\infty<t<0$.
By considering higher derivatives in time of these equations and  Proposition \ref{bootstrapintime}, one can obtain Theorem \ref{mabssmooth} using the regularity theory of the stationary Stokes system together with bootstrap arguments. $\Box$

 \setcounter{equation}{0}
\section{Proof of Theorem \ref{hsmildsol}}



 \begin{lemma}\label{smoothms} Assume that a bounded function $u$ satisfies conditions (\ref{log}), (\ref{p2}), (\ref{incompressibility in hs}), and (\ref{momentumpres}) of Theorem \ref{hsmildsol}. Then $\nabla u\in L_\infty(Q_-^+)$. The function $u$ is infinitely smooth in spatial variables in upper half space $x_3>0$. \end{lemma}
 \textsl{Proof} Let $z_0=(x'_0,0,t_0)$,  $Q_+(z_0,R)=B_+(x_0,R)\times ]t_0-R^2,t_0[$, and $R_k=R-R\sum\limits_{i=1}^k2^{-i-1}$ for $k=1,2,...$, and $R_0=R$.
Let $\varphi\in C^\infty_0(B(x_0,R)\times ]t_0-R^2,t_0+R^2[)$
and let $v:=u\varphi$ and  $R=1$. Then
$$\partial_tv-\Delta v=f_1+f_2,$$
where
$$f_1:=-\varphi \nabla p^2+u(\partial_t \varphi+\Delta\varphi)+ u\cdot\nabla\varphi u+(p^1-[p^1]_{B(x_0,1)})\nabla \varphi$$
and
$$f_2:=-{\rm div}(2u\otimes \nabla\varphi +\varphi(u\otimes u+(p^1-[p^1]_{B(x_0,1)})\mathbb I)).$$
Moreover, $v=0$ satisfies $\partial'Q_+(z_0,1)$. We can split $v$ into two parts $v=v^1+v^2$ so that
$$\partial_tv^1-\Delta v^1=f_1$$
and
$v^1=0$ on $\partial'Q_+(z_0,1)$. By our assumptions,
$$\|f_1\|_{s,\infty,Q_+(x_0,1)}\leq c(s)$$
for any $1<s<\infty$. Therefore we can claim that
$$|\nabla v^1|\leq c$$
on $Q_+(z_0,3/4)$ with a constant independent of $z_0$.

Notice that by our assumptions we may write $f_2= {\rm div}(F_{2})$, where $$\|F_2\|_{s,\infty,Q_{+}(x_{0},1)}\leqslant c(s).$$
So using boundary regularity theory for the heat equation, we can say that
$$\|\nabla v^2\|_{s, Q(z_0,3/4)}\leq c.$$
Then we can see that, since
$$\Delta  p^1=-{\rm div} (u\cdot\nabla u)$$
with Neumann boundary condition on the flat part of the boundary,
$$\|\nabla p^1\|_{s,Q_+(z_0, (3/4+5/8)/2)}\leq c.$$
This means that
$$\|\partial_t u\|_{s,Q_+(5/8)}+\|\nabla^2u\|_{s,Q_+(5/8)}\leq c$$
and thus by the parabolic imbedding theorem (for large enough $s$)

$$|\nabla u|\leq c$$
on $Q_+(z_0,1/2)$ with a constant independent of $z_0$.
So, we have boundedness near the boundary. To get the interior estimate, we can use the same arguments that used for mild bounded ancient solutions in the whole space.

In fact, we have even more
$\nabla u$ is continuous up the boundary and $$\|\partial_tu\|_{s,Q(z_0,1)\cap Q^+_-}+\|\nabla^2u\|_{s,Q(z_0,1)\cap Q^+_-}\leq c(s)$$ for any $s>1$.
Lemma \ref{smoothms} is proven. $\Box$


We wish to show that $u$ has the following properties. For any $A<0$,
$$u=u^1+u^2,$$
where
$$ u^1(x,t)=\int\limits_{\mathbb R^3_+}G(x,y,t-A)u(y,A)dy$$  and $$u^{2}_i(x,t)=\int\limits^t_A\int\limits_{\mathbb R^3_+}K_{ijm}(x,y,t-\tau)u_j(y,\tau)u_{m}(y,\tau)dy d\tau
$$ in $Q_{A}$.

Let us go back to Lemma \ref{itegbyparts} and its proof.  $F$  and its approximations $F^m$ are from that lemma an its proof.  Solonnikov showed in \cite{Sol2003} that $v^m$  given by the formula

$$v^m(x,t)=\int\limits_{\mathbb R^3_+}G(x,y,t-A)u(y,A)dy+\int\limits^t_A\int\limits_{\mathbb R^3_+}K(x,y,t-s)F^m(y,s)dy ds=$$$$=\int\limits_{\mathbb R^3_+}G(x,y,t)u(y,A)dy+\int^t_A\int\limits_{\mathbb R^3_+}G_{ij}(x,y,t-s)f^m_j(y,s)dy$$
 satisfies the identity:
$$\int\limits_{Q_A^+}v^m\cdot\nabla q dz=0$$
for any $q\in C^\infty_0(Q_A)$, where $Q_A^+:=\mathbb R^3_+\times ]A,0[$  and $Q_A:=\mathbb R^3\times ]A,0[$,
$$
    \int\limits_{Q^+_A}v^m\cdot(\partial_t\varphi+\Delta \varphi)dx dt+\int\limits_{\mathbb R^3_+}u(x,A)\cdot \varphi(x,A)dx=  -\int\limits_{Q^+_A}f^m\cdot \varphi  dx dt
$$
for any divergence free functions $\varphi\in C^\infty_{0}(Q_-)$ with $\varphi(x',0,t)=0$ for any $x'\in \mathbb R^2$ and for any $t<0$.  Since
$$\int\limits_{\mathbb R^3_+}|G(x,y,t)|dy\leq c,$$
we can use boundedness and pass to the limit as $m\to\infty$.

As a result, we have
$$\int\limits_{Q_A^+}v\cdot\nabla q dz=0$$
for any $q\in C^\infty_0(Q_A)$, where  $Q_A:=\mathbb R^3\times ]A,0[$. Here, $v$ is defined as $u_{1}+u_{2}$.  Furthermore, for any divergence free functions $\varphi\in C^\infty_{0}(Q_-)$ with $\varphi(x',0,t)=0$ for any $x'\in \mathbb R^2$ and for any $t<0$:
$$
    \int\limits_{Q^+_A}v\cdot(\partial_t\varphi+\Delta \varphi)dx dt+\int\limits_{\mathbb R^3_+}u(x,A)\cdot \varphi(x,A)dx=  -\int\limits_{Q^+_A}f\cdot \varphi  dx dt.
$$

Now assume that $F=u\otimes u$.
From the the previous pages, it is clear that $u$ is continuous in the completion of $Q_+(R):=B_+(R)\times ]-R^2,0[$. Then using cut-off functions in time, we can show that $u$ satisfies the same identity as $v$. And thus letting
$w=u-v$, we get
$$\int\limits_{Q_A^+}w\cdot\nabla q dz=0$$
for any $q\in C^\infty_0(Q_A)$,
$$
    \int\limits_{Q^+_A}w\cdot(\partial_t\varphi+\Delta \varphi)dx dt=  0
$$
for any divergence free functions $\varphi\in C^\infty_{0}(Q_-)$ with $\varphi(x',0,t)=0$ for any $x'\in \mathbb R^2$ and for any $t<0$.
If we extend $w$ by zero for $t\leq A$, we find
$$\int\limits_{Q_-^+}w\cdot\nabla q dz=0$$
for any $q\in C^\infty_0(Q_-)$,
$$
    \int\limits_{Q^+_-}w\cdot(\partial_t\varphi+\Delta \varphi)dx dt=  0
$$
for any divergence free functions $\varphi\in C^\infty_{0}(Q_-)$ with $\varphi(x',0,t)=0$ for any $x'\in \mathbb R^2$ and for any $t<0$.

By the Liouville theorem (see \cite{JiaSS2012} and \cite{JiaSS2013}), $w=w(x_3,t)$. We need to show that $w\equiv0$ in $Q_A^+$.
To this end, it is sufficient to show that for $x^{'}\in\mathbb{R}^2$ and $t\in ]A,0[$, one has $\nabla q^2(x',x_{3},t)\to 0$ as $x_3\to \infty$. Here, $q^2$ is the pressure for $v$ so that
$$\partial_tv-\Delta v+\nabla q^2=f$$
in $Q_A^+$. If split $v=v^1+v^2$ so that $v^i$ corresponds to the Green function $G^i$. Then, clearly,
$$\partial_tv^1-\Delta v^1=f$$
in $Q_A^+$,
$$v^1(x',0,t)=0$$
for $x'\in\mathbb R^2$ and $A<t<0$, and
$$v^1(\cdot,0)=u_0(\cdot).$$
Thus,
$$\nabla q^2=\Delta v^2-\partial_tv^2.$$
On the other hand, we have $v^2=v^{2,1}+v^{2,2}$, where
$$v^{2,1}(x,t)=\int\limits_{\mathbb R^3}G^2(x-y,t-A)u_0(y)dy$$ and
$$v^{2,2}(x,t)=\int\limits_A^t\int\limits_{\mathbb R^3}G^2(x-y,t-s)f(y,s)dy ds.$$
Using reasoning from Proposition \ref{propertytimeder}, obtain
$$|\partial_tv^{2,1}(x,t)|\leq  \frac c{t-A}\int\limits_{\mathbb R^3_+}\frac 1{(|x'-y'|^2+x_3^2+y_3^2+t-A)^\frac 32}e^{-\frac {cy^2_3}{t-A}}dy'dy_3\leq $$
$$\leq  \frac c{t-A}\int\limits_0^\infty \frac 1{(x_3^2+y_3^2+t-A)^\frac 12}e^{-\frac {cy^2_3}{t-A}}dy_3\leq $$
$$\leq   c(t-A)^{-\frac 12}((x_3^2+t-A)^{-\frac 12}\to0$$
as $x_3\to\infty$.
Next, since $G^{2}(x,y,0)=0$, then
$$|\partial_t v^{2,2}(x,t)|\leq c\int\limits_A^t(t-s)^{-\frac 12}((x_3^2+t-s)^{-\frac 12}ds\to0$$
as $x_3\to\infty$.

Regarding $v^{2,2}$, we have
$$|\nabla ^2v^{2,1}(x,t)|\leq \frac c{t-A+x^2_3}\int\limits_{\mathbb R^3_+} \frac 1{(|x'-y'|^2+x_3^2+y_3^2+t-A)^\frac 32}e^{-\frac {cy^2_3}{t-A}}dy'dy_3\leq $$
$$\leq  c(t-A+x^2_3)^{-\frac 12}((x_3^2+t-A)^{-\frac 12}= c(t-A+x^2_3)^{-1}\to0$$
as $x_3\to\infty$.
Next, $$|\nabla ^2v^{2,2}(x,t)|\leq c\int\limits^t_A(t-s+x^2_3)^{-1}ds\to0$$
as $x_3\to\infty$. So, we have the required decay for $\nabla q^2$ and thus we have for all $A<0$ the following integral representation:
$$u(x,t):=\int\limits_{\mathbb R^3}G(x,y,t-A)u(y,A)dy+$$$$+\int\limits_A^t\int\limits_{\mathbb R^3}G(x,y,t-s){\rm div} (u\otimes u-p^1\mathbb I)(y,s)dy ds=$$$$=\int\limits_{\mathbb R^3}G(x,y,t-A)u(y,A)dy+\int\limits^t_A\int\limits_{\mathbb R^3_+}K(x,y,t-\tau)u(y,\tau)\otimes u(y,\tau)dy d\tau$$
for all $t>A$.

Now, our aim is to prove the inverse statement, i.e., we assume that bounded divergence free function satisfied the latter identity for any $A<0$. Introducing $F=u\otimes u$ and tensor $H=F+p^1\mathbb I$ and using approximations of $u$,
we can show that $u$ is a distributional solution to the Navier-Stokes equations in $Q^+_-$ and belongs to the space $W^1_\infty(Q^+_-)$. This can be done in the same way as in \cite{SerSve2013} (the most difficult part of that paper).
We then can introduce the pressure $p^2$
 so that $$\partial_tu-\Delta u+\nabla p^2=-{\rm div}H.$$ Splitting $u=u^1+u^2$ and repeating the aforesaid arguments, we can show that $\nabla p^2$ satisfies all requirements in the definition of bounded mild ancient solutions.

\setcounter{equation}{0}
\section{Proof of Proposition \ref{equivalence in space}}
Assume first that $u$ satisfies the conditions of Proposition \ref{equivalence in space}, i.e.,
there exists a pressure field $p\in L_\infty(-\infty,0;BMO(\mathbb R^3))$ such that
\begin{equation}\label{momentumwithpressure}
    \int\limits_{Q_-}\Big(u\cdot(\partial_t\varphi+\Delta \varphi)+u\otimes u:\nabla\varphi\Big)dz=-\int\limits_{Q_-}p\,{\rm div}\,\varphi dz.
    \end{equation}

    Our aim is to show that, for any $A<0$, the function satisfies the integral identity (\ref{mildinwholespace}).
    First,  let us notice that the presure $p$ (up to a bounded function of time $t$) is formally represented as follows:
 $$p(x,t)=-\frac 13 |u(x,t)|^2+\frac 1{4\pi}\int\limits_{\mathbb R^3} \nabla^2_y\Big(\frac  1{|x-y|}\Big):u(y,t)\otimes u(y,t)dy.$$




We know that mild bounded ancient solutions are infinitely smooth and all partial derivatives, apart from
derivatives in time of the pressure, are bounded. The derivatives $\partial^k_tp$, $k=0,1,...$, belong to
$L_\infty(BMO)$, see \cite{Ser2014}. So, we re-write the Navier-Stokes equations in the following way:
$$\partial_tu-\Delta u=f,\qquad {\rm div}u=0$$
in $Q_-$, where $f:=-{\rm div}\, u\otimes u-\nabla p$. We know that $f$ is infinitely smooth and all its derivative are bounded.  Then, by Tychonoff's uniqueness theorem, 
$$u(x,t):=\int\limits_{\mathbb R^3}\Gamma(x-y,t)u(y,A)dy+\int\limits_A^t\int\limits_{\mathbb R^3}\Gamma(x-y,t-s)f(y,s)dy ds$$
for $t> A$ and for all $A<0$. It remains to show
that
$$\int\limits_A^t\int\limits_{\mathbb R^3}\Delta_y\Phi(x-y,t-s)f_i(y,s)dyds=
\int\limits_A^t\int\limits_{\mathbb R^3}\Gamma(x-y,t-s)f_i(y,s)dyds=$$$$=\int\limits^t_A \int\limits_{\mathbb R^3} K_{ijm}(x,y,t-\tau)u_j(y,\tau)u_m(y,\tau)dy d\tau$$
for any $A<0$.
To this end, we introduce as notation $p_F$ which is the BMO-solution to
to equation $\Delta p=-{\rm div }{\rm div}F $ in $\mathbb R^3$ with $[p_F]_B=0$. We deduce the required identity from the following lemma.
\begin{lemma}\label{identi}
Let $F$ be a bounded smooth function in $\mathbb R^3$ having all derivatives bounded there. Then, for any positive $s$,
$$\int\limits_{\mathbb R^3}\Delta_y\Phi(x-y,s)f_i(y)dy=
\int\limits_{\mathbb R^3}\Gamma(x-y,s)f_i(y)dy=$$$$=\int\limits_{\mathbb R^3} K_{ijm}(x,y,s)F_{jm}(y)dy, $$
where $f=-{\rm div} F-\nabla p_F$.
\end{lemma}
 The proof can be done with the help of suitable approximation of $F$ and the following estimates:
$$|\nabla^k_x\Phi(x,t)|\leq \frac {c(k)}{(t+|x|^2)^\frac {1+k}2}$$ and
$$|\nabla^k\Gamma(y,1)|\leq \frac {c(k)}{(1+|y|^2)^\frac {3+k}2}e^{-\frac {|y|^2}8}. $$
 The first one is due to Solonnikov, see \cite{Sol1964}), and the second one is well known.

Inverse statement of Proposition \ref{equivalence in space} can be easily deduced  from the above lemma and suitable approximations of $u$. This completes the proof of Proposition \ref{equivalence in space}.

\setcounter{equation}{0}
\section{Appendix }
 \textbf{\textsl{Proof of Lemma \ref{Uniforminteggeneralise}}}
  For $x'$ in $\mathbb{R}^{2}$, denote by $Q_{\mathbb{R}^{2}}(x^{'},R)$  a cube in $\mathbb{R}^{2}$ with side lengths $2R$ centred at $x'$. Define the   space cylinders and the space-time cylinders:
$$C(x',R,m):= Q_{\mathbb{R}^{2}}(x^{'},R)\times ]m,m+1[,$$
$$C(x',R,m,A,t):= C(x',R,m) \times ]A,t[.$$
 After a decomposition of the domain,  consider the following integrals separately ($R=1,2,3\ldots$ and $m=0,1,2\ldots$):
\begin{equation}\label{nearsingularity1}
I(m,1,x,t):= \int\limits_{C(x',1,m,A,t)}|\nabla_{y}G^{2}(x,y,t-\tau)f(y,\tau)| dyd\tau
\end{equation}
\begin{equation}\label{farsingularity1}
I(m,R,x,t):= \int\limits_{C(x',R+1,m,t,A)\setminus C(x',R,m,t,A)} |\nabla_{y}G^{2}(x,y,t-\tau)f(y,\tau)| dyd\tau.
\end{equation}
First consider $I(0,1,x,t)$.
Let
$$ J(x,y,t-\tau)=(|x-y^*|^2+t-\tau)^{-\frac {3}2}\exp{\Big(-\frac {cy^2_3}{t-\tau}\Big)}.$$
 From the Solonnikov estimates (\ref{GreenFunctEst2}):
$$I(0,1,x,t) \leqslant c\int\limits_{C(x',1,0,A,t)}(t-\tau)^{-\frac{1}{2}}J(x,y,t-\tau)|f(y,\tau)| dyd\tau.$$
Then by the H\"older inequality we have
$$|I(0,1,x,t)|^{l'}\leqslant c\|f\|_{L_{s,l,unif}(Q_{A}^{+})}^{l'}\int\limits^t_{A}\Big(\int\limits_{C(x',1,0)}(t-\tau)^{-\frac{s'}{2}}|J(x,y,t-\tau)|^{s'} dy\Big)^{\frac{l'}{s'}}d\tau.$$
We get after a change of variables
$$\int\limits^t_{A}\Big(\int\limits_{C(x',1,0)}(t-\tau)^{-\frac{s'}{2}}|J(x,y,t-\tau)|^{s'} dy\Big)^{\frac{l'}{s'}}d\tau$$$$ \leqslant \int\limits^t_A(t-\tau)^{-2l'+
\frac{3l'}{2s'}}\Big(\int\limits_{\mathbb{R}^{3}_{+}}(|z|^{2}+1)^{\frac{3s'}{2}}\exp{(- {cs'z^2_3})} dz\Big)^{\frac{l'}{s'}}d\tau\leqslant$$$$\leqslant c(s,l)\int\limits^{-A}_{0}\lambda^{-2l'+\frac{3l'}{2s'}}d\lambda.
$$
This quantity is finite if and only if (\ref{Uniformintegcondition}) holds.

For $I(m,1,x,t)$, with $m\geqslant 1$, the H\"older inequality gives
$$|I(m,1,x,t)|^{l'}\leqslant c\|f\|_{L_{s,l,unif}(Q_{A}^{+})}^{l'}\int\limits^t_{A}\Big(\int\limits_{C(x^{'},1,m)}(t-\tau)^{-\frac{s'}{2}}|J(x,y,t-\tau)|^{s'} dy\big)^{\frac{l'}{s'}}d\tau.$$
For the second factor we have $$\int\limits^t_{A}\Big(\int\limits_{C(x',1,m)}(t-\tau)^{-\frac{s'}{2}}|J(x,y,t-\tau)|^{s'} dy\big)^{\frac{l'}{s'}}d\tau\leqslant$$$$
\leqslant\int\limits^t_A\Big(\int\limits_{Q_{\mathbb{R}^{2}}(0,1)}\frac{\exp{\Big(-\frac {cm^2s'}{t-\tau}\Big)}}{{(t-\tau)}^{\frac{s'}{2}}}(|y'|^{2}+t-\tau)^{\frac{-3s'}{2}} dy'\Big)^{\frac{l'}{s'}}d\tau\leqslant \frac{C(s,A)}{m^{2}}.$$
Here for the final line the following fact  is used (for $\alpha>0$):
$$\sup_{x>0}x^{\alpha}\exp(-x)\leqslant C(\alpha).$$

Now, consider $I(0,R,x,t)$. Initially using same arguments as for $I(0,1,x,t)$, we have
$$|I(0,R,x,t)|^{l'}\leqslant \|f\|_{s,l,C(x^{'},R+1,0,A,t)\setminus C(x', R,0,A,t)}^{l'}\times$$$$\times\int\limits^t_A\Big(  \int\limits_{C(x^{'},R+1,0)\setminus C(x^{'},R,0)}(t-\tau)^{-\frac{s'}{2}}J(x,y,t-\tau)^{s'} dy\,\Big)^{\frac{l'}{s'}}d\tau\leqslant $$$$\leqslant c(s,l) R^{\frac{l'}{s}}\|f\|^{l'}_{L_{s,l,unif}(Q_{A}^{+})}\times$$$$\times
\int\limits^t_A\Big(  \int\limits_{C(x',R+1,0)\setminus C(x',R,0)}(t-\tau)^{-\frac{s'}{2}}J(x,y,t-\tau)^{s'} dy\,\big)^{\frac{l'}{s'}}d\tau\leqslant$$$$\leqslant c(s,l) R^{\frac{l'}{s}}\|f\|^{l'}_{L_{s,l,unif}(Q_{A}^{+})}\times$$$$\times
\int\limits^t_{A}(t-\tau)^{\frac{-l'}{2}
}\Big(\int\limits_{Q_{\mathbb{R}^{2}}(x',R+1)\setminus Q_{\mathbb{R}^{2}}(x',R)}|x'-y'|^{-3s'}dy'\Big)^{\frac{l'}{s'}}d\tau\leqslant$$$$\leqslant c(s,l)R^{-3l'+\frac{2l'}{s'}+\frac{l'}{s}}\|f\|^{l'}_{L_{s,l,unif}(Q_{A}^{+})}\int\limits^{-A}_{0}\lambda^{\frac{-l'}{2}
}
d\lambda. $$
  By (\ref{Uniformintegcondition}), $l>2$ and thus $l'<2$. So, the last factor is finite.
 Hence,
 $$ |I(0,R,x,t)|\leqslant c(A,s,l)R^{-1-\frac{1}{s}}\|f\|_{L_{s,l,unif}(Q_{A}^{+})}.$$
 Similar arguments to before give, for $m\geqslant 1$,
 $$|I(m,R,x,t)|\leqslant\frac{c(A,s,l)}{m^{2}} R^{-1-\frac{1}{s}}\|f\|_{L_{s,l,unif}(Q_{A}^{+})}.$$
 Summing over $m$ and $R$ we then conclude. $\Box$

\begin{lemma}\label{singularintrgal}
Let $n\geqslant 3$. Denote $\Delta:= \{(x,y)\in\mathbb{R}^{n}\times\mathbb{R}^{n}: x\neq y\}.$
Suppose the measurable kernel $K:\mathbb{R}^{n}\times\mathbb{R}^{n}\setminus\Delta\rightarrow \mathbb{R}$ is such that there exists $M>0$ with:
\begin{equation}\label{boundedkernelgeneral}
 |K(x,y)|\leqslant\frac{M}{|x-y|^{n}}.
\end{equation}
Define the truncation (on $L_{p}(\mathbb{R}^{n})$, $1<p<\infty$):
\begin{equation}\label{truncation}
T_{\epsilon}(f)(x):= \int\limits_{|x-y|\geq\varepsilon}K(x,y)f(y)dy.
\end{equation}
Suppose, for this kernel,  there exists bounded linear operator $T:L_{p}(\mathbb{R}^{n})\rightarrow L_{p}(\mathbb{R}^{n})$ ($1<p<\infty$) such that:
that for $f\in L_{p}(\mathbb{R}^{n})$ ($1<p<\infty$):
\begin{equation}\label{convergtruncation}
\|T_{\epsilon}(f)-T(f)\|_{L_{p}(\mathbb{R}^{n})}\rightarrow 0,
\end{equation}
\begin{equation}\label{Lpbounded}
\|T(f)\|_{L_{p}(\mathbb{R}^{n})}\leqslant c(K,n)\|f\|_{L_{p}(\mathbb{R}^{n})}.
\end{equation}
Furthermore for $f\in L_{\infty}(\mathbb{R}^{n})$ compactly supported:
\begin{equation}\label{LinftyBMO}
\|T(f)\|_{BMO(\mathbb{R}^{n})}\leqslant c(K,n)\|f\|_{L_{\infty}(\mathbb{R}^{n})}.
\end{equation}
Here, $C(K,n)$ means that the constant depends on the properties of the Kernel (e.g some smoothness of the kernel) and the dimension of the space.\\
Consider an unbounded domain $\Omega\subset\mathbb{R}^{n}$ that is contained between two $n-1$ dimensional parallel hyperplanes (denoted $\Pi_{1}$ and $\Pi_{2}$ respectively) a finite distance $2L$ apart.
Take $g$ to be a compactly supported function in $L_{p}(\mathbb{R}^{n})$ (for $1<p<\infty$) such that $g$ is non-zero and bounded outside of $\Omega$.\\
Then it follows that:
\begin{equation}\label{decomp}
 Tg(x)= h_{1}(x)+h_{2}(x).
 \end{equation}
 Here,
 \begin{equation}\label{BMOpart}
 \|h_{1}\|_{BMO(\mathbb{R}^{n})}\leqslant c(K,n)\|g\|_{L_{\infty}(\mathbb{R}^{n}\setminus\Omega)}
 \end{equation}
 and
 \begin{equation}\label{Luniformpart}
 \|h_{2}\|_{L_{p,unif}(\mathbb{R}^{n})}\leqslant c(K,M,n,p,L)\|g\|_{L_{p,unif}(\mathbb{R}^{n})}.
 \end{equation}
\end{lemma}

\textsl{Proof of Lemma \ref{singularintrgal}}\
 For $x'$ in $\mathbb{R}^{n-1}$, denote by $Q_{\mathbb{R}^{n-1}}(x^{'},R)$  a cube in $\mathbb{R}^{n-1}$ with side lengths $2R$ centred at $x'$. \\First one shows that, without loss of generality, it is sufficient to reduce to the case where:
\begin{equation}\label{rotatedplane1}
\Pi_{1}= \{(x',L):x'\in\mathbb{R}^{n-1}\},
\end{equation}
\begin{equation}\label{rotatedplane2}
\Pi_{2}= \{(x',-L): x'\in\mathbb{R}^{n-1}\}.
\end{equation}
Let $A:\mathbb{R}^{n}\rightarrow\mathbb{R}^{n}$ be a rotation of $\mathbb{R}^{n}$.
It can be inferred that:
\begin{equation}\label{truncationrotation}
 T_{\epsilon}(g)(A(x))= \int\limits_{|x-y|\geq\varepsilon}K(A(x),A(y))g(A(y))dy.
 \end{equation}
 If one lets
 \begin{equation}\label{rotatedkernel}
 \bar{K}(x,y):= K(A(x),A(y)),
 \end{equation}
clearly $\bar{K}$ satisfies  (\ref{boundedkernelgeneral}).
Define the  truncation operator (for $f\in L_{p}(\mathbb{R}^{n})$, $1<p<\infty$):
\begin{equation}\label{rotatedtruncation}
S_{\epsilon}(f)(x):= \int\limits_{|x-y|\geqslant\epsilon}\bar{K}(x,y)f(y) dy.
\end{equation}
By rotation invariance of $L_{p}(\mathbb{R}^{n})$ and $BMO(\mathbb{R}^{n})$, it is inherited from $T_{\epsilon}$ and $T$ that there exists a bounded linear operator $S:L_{p}(\mathbb{R}^{n})\rightarrow L_{p}(\mathbb{R}^{n})$ (where $1<p<\infty$) such that (\ref{convergtruncation})-(\ref{LinftyBMO}) hold.
 Since the space $L_{p,unif}(\mathbb{R}^{n})$ is rotation invariant, rotations of  $\Omega$ can be considered without loss of generality.\\
 Fix $z_{1}\in\mathbb{R}^{n}$. It can be inferred that:
 \begin{equation}\label{truncationtranslation}
  T_{\epsilon}(g)(x-z_{1})= \int\limits_{|x-y|\geq\varepsilon}K(x-z_{1},y-z_{1})g(y-z_{1})dy.
 \end{equation}
 Let
 \begin{equation}\label{translatedkernel}
 \bar{K}(x,y):= K(x-z_{1},y-z_{1}).
 \end{equation}
 Using the spaces $BMO(\mathbb{R}^{n})$, $L_{p}(\mathbb{R}^{n})$ and $L_{p,unif}(\mathbb{R}^{n})$ are translation invariant,  one can use the aforementioned arguments to show that translations of  $\Omega$ can be considered without loss of generality.\\
 From now on take $\Pi_{1}$ as in (\ref{rotatedplane1}) and $\Pi_{2}$ as in (\ref{rotatedplane2}).\\
   Now decompose $g$:
 $$g_{1}(x',x_{n})= (1-\chi_{]-L,L[}(x_{n}))g(x',x_{n}),$$
 $$g_{2}(x',x_{n})= \chi_{]-L,L[}(x_{n})g(x',x_{n}).$$
 By (\ref{LinftyBMO}), get that $h_{1}(x):= T(g_{1})(x)$ satisfies (\ref{BMOpart}). It remains to show $h_{2}(x):= T(g_{2}(x))$ satisfies (\ref{Luniformpart}).
 By an identical argument to that showing translations of $\Omega$ are permissible, it is sufficient to prove:
 \begin{equation}\label{Luniformsimplification}
 \sup_{x_{n}\in\mathbb{R}}\|h_{2}\|_{L_{p}(B((0,x_{n}),1))}\leqslant c(K,M,n,p,L)\|g\|_{L_{p,unif}(\mathbb{R}^{n})}.
 \end{equation}
 Let us write $g_{2}:= g^+_{2}+g^-_{2}$,
 where
 \begin{equation}\label{g2decomp}
 g^-_{2}(x',x_{3})=\chi_{Q_{\mathbb{R}^{n-1}}(0,2)}(x')g_{2}(x',x_{n}).
 \end{equation}
 Further to this write $h_{2}:= h^+_{2}+h^-_{ 2}$, where
 \begin{equation}\label{h2decomp}
 h^-_{2}= T(g^-_{2}).
 \end{equation}
 By (\ref{Lpbounded}), $h^-_{2}$
 satisfies the estimate (\ref{Luniformsimplification}) in place of $h_{2}$. It remains to show the same for $h^+_{2}$.
 It can be shown (for $z$ in $B((0,x_{n}),1)$),
  $$|h^+_{2}(z)|\leqslant \int\limits^L_{-L}\int\limits_{y'\in\mathbb{R}^{n-1}\setminus Q_{\mathbb{R}^{n-1}}(0,2)}\frac{M}{|z'-y'|^{n}}|g(y^{'},y_{n})|dy'dy_{n}\leqslant$$$$\leqslant c(n,M) \int\limits^L_{-L}\int\limits_{y'\in\mathbb{R}^{n-1}\setminus Q_{\mathbb{R}^{n-1}}(0,2)}\frac{1}{|y'|^{n}}|g(y',y_{n})|dy'dy_{n}
 \leqslant$$$$\leqslant c(n,M)\sum_{N=2}^{\infty}\int\limits^L_{-L}\int\limits_{Q_{\mathbb{R}^{n-1}}(0,N+1)\setminus Q_{\mathbb{R}^{n-1}}(0,N)}\frac{1}{|y'|^{n}}|g(y',y_{n})|dy'dy_{n}.$$
 The domain $Q_{\mathbb{R}^{n-1}}(0,N+1)\setminus Q_{\mathbb{R}^{n-1}}(0,N)\times ]-L,L[$, can be seen to be covered by $c(n)N^{n-2}\times\lceil L \rceil$ unit cylinders. Here $\lceil L \rceil$ is the smallest integer greater than $L$. Hence, by H\"older's inequality:

 $$\int\limits^L_{-L}\int\limits_{Q_{\mathbb{R}^{n-1}}(0,N+1)\setminus Q_{\mathbb{R}^{n-1}}(0,N)}\frac{1}{|y'|^{n}}|g(y',y_{n})|dy'dy_{n} \leqslant
$$$$\leqslant
  c(n,L,p)N^{\frac{n-2}{p}}\|g\|_{L_{p,unif}(\mathbb{R}^{n})}\Big(\int\limits_{Q_{\mathbb{R}^{n-1}}(0,N+1)\setminus Q_{\mathbb{R}^{n-1}}(0,N)}\frac{1}{|y'|^{np'}}dy'dy_{n}\Big)^{\frac{1}{p'}}.
 $$
One can  estimate for the second factor and get the bound:
$$\frac{c(n,p)}{N^{n-\frac{n}{p'}+\frac{1}{p'}}}.$$
Thus,
$$\int\limits^L_{-L}\int\limits_{Q_{\mathbb{R}^{n-1}}(0,N+1)\setminus Q_{\mathbb{R}^{n-1}}(0,N)}\frac{1}{|y'|^{n}}|g(y',y_{n})|dy'dy_{n} \leqslant$$$$\leqslant c(n,L,p)\|g\|_{L_{p,unif}(\mathbb{R}^{n})}N^{-(1+\frac{1}{p})}.$$
So it is obtained that:
$$ |h_{2}^{+}(z)|\leqslant c(n,M,L,p)\|g\|_{L_{p,unif}(\mathbb{R}^{n})}\sum_{N=2}^{\infty} N^{-(1+\frac{1}{p})}.$$
From here all conclusions follow immediately. $\Box$




\end{document}